\documentclass[12pt]{amsart}
\usepackage{amsmath, amssymb, latexsym}

\usepackage[dvips]{graphicx}
\usepackage[dvips]{color}

\title[Equivalence of domains ...]{Equivalence of domains with isomorphic 
semigroups of endomorphisms}

\author[S.~Merenkov]{Sergei Merenkov}
\address{Department of Mathematics\\ Purdue University\\ West Lafayette, 
Indiana 47907}
\email{smerenko@math.purdue.edu}
\thanks{Research supported by NSF, DMS 0072197}

\newcommand\C{{\mathbb C}}
\newcommand\R{{\mathbb R}}

\newtheorem{theorem}{Theorem}

\newtheorem{lemma}{Lemma}
\newtheorem{remark}{Remark}

\begin{document}


\abstract{For two bounded domains $\Omega_1,\ \Omega_2$ in $\C$ whose 
semigroups of analytic endomorphisms $E(\Omega_1), \ E(\Omega_2)$ are 
isomorphic with an isomorphism 
$\varphi:\ E(\Omega_1)\rightarrow E(\Omega_2)$, Eremenko proved in 1993 that
there exists a conformal 
or anticonformal map $\psi:\ \Omega_1\rightarrow \Omega_2$ such that
$
\varphi f=\psi\circ f\circ \psi^{-1}, 
$
for all  $f\in E(\Omega_1)$.

In the present paper we prove an analogue of this result for the case of 
bounded domains in $\C^n$.}
\endabstract

\maketitle

\section{Introduction}\label{S:Intro}

A classical theorem of L. Bers says that every $\C$-algebra isomorphism 
$H(A)\rightarrow H(B)$ of algebras of holomorphic functions in domains $A$ 
and $B$ in the complex plane has the form $f\mapsto f\circ\theta$, where 
$\theta:\ B\rightarrow A$ is a conformal isomorphism, or $f\mapsto 
\overline{f}\circ\theta$ with anticonformal $\theta$. In particular, the 
algebras $H(A)$ and $H(B)$ are isomorphic if and only if the domains $A$ 
and $B$ are conformally equivalent. H. Iss'sa \cite{hI66} obtained a similar 
theorem for fields of meromorphic functions on Stein spaces. A good 
reference for these results is 
\cite{mH68}.

Likewise, a question of recovering a topological space from the algebraic 
structure of its semigroup of continuous self-maps has been extensively
studied \cite{kM75}.

In 1990, L. Rubel asked whether similar results hold for semigroups (under 
composition) $E(D)$ of holomorphic endomorphisms of a domain $D$. 
A. Hinkkanen constructed 
examples \cite{aH92} which show that even non-homeomorphic domains in 
$\C$ can have isomorphic semigroups of endomorphisms. The reason is 
that the semigroup of endomorphisms of a domain can be too small to 
characterize this domain. 

However, in 1993, A. Eremenko \cite{aE93} proved that for two Riemann 
surfaces $D_1$, $D_2$, which admit bounded nonconstant holomorphic 
functions, and such that the semigroups of analytic endomorphisms $E(D_1)$ 
and $E(D_2)$ are isomorphic with an isomorphism 
$\varphi: \ E(D_1)\rightarrow E(D_2)$, there exists a conformal or 
anticonformal map $\psi:\ D_1\rightarrow D_2$ such that 
$\varphi f=\psi\circ f\circ\psi^{-1}$, for all $f\in E(D_1)$. 
In the present paper we investigate the analogue of this result 
for the case of bounded domains in $\C^n$. The theorems of Bers 
and Iss'sa, mentioned above, do not extend to arbitrary domains in 
$\C^n$.

For a bounded domain $\Omega$ in $\C^n$ we denote by $E(\Omega)$ the 
semigroup of analytic endomorphisms of $\Omega$ under composition. In what 
follows, we say that a map is {\it{(anti-) biholomorphic}}, if it is 
biholomorphic 
or antibiholomorphic. We prove the following theorem.

\begin{theorem}\label{T:Mt}
Let $\Omega_1,\ \Omega_2$ be bounded domains in $\C^n,\ \C^m$ 
respectively, and suppose that there exists 
$\varphi:\ E(\Omega_1)\rightarrow E(\Omega_2)$,  an isomorphism of 
semigroups. Then $n=m$ and there exists an (anti-) 
biholomorphic map $\psi:\ \Omega_1\rightarrow\Omega_2$ such that 
\begin{equation}\label{E:C}
\varphi f=\psi\circ f\circ \psi^{-1}, \ \ \text{for all}\ f\in E(\Omega_1). 
\end{equation} 
\end{theorem}

The existence of a homeomorphism $\psi$ satisfying (\ref{E:C}) follows from 
simple general considerations (Section~\ref{S:Top}). The hard part is proving 
that 
$\psi$ is (anti-) biholomorphic. In dimension 1 this is done by 
linearization of holomorphic germs of $f\in E(\Omega)$ near an attracting 
fixed point. In several dimensions such linearization theory exists 
(\cite{vA88}, pp. 192--194), but it is too complicated (many germs with 
an attracting fixed point are non-linearizable, even formally). 
In Sections~\ref{S:Loc},~\ref{S:Ext} we show how to localize the problem. 
In Sections~\ref{S:Sys},~\ref{S:Sim}  
we describe, using only the semigroup structure, a large enough class of 
linearizable germs. Linearization of these germs permits us to reduce the 
problem to a matrix functional equation, which is solved in 
Section~\ref{S:Sol}. In Section~\ref{S:Pr} we complete the proof that
$\psi$ is (anti-) biholomorphic.

Theorem~\ref{T:Mt} can be slightly generalized, namely one may assume
that $\varphi$ is an epimorphism. In Section~\ref{S:Pr2} we prove the
following theorem.

\begin{theorem}\label{T:Gt}
If $\varphi:\ E(\Omega_1)\rightarrow E(\Omega_2)$ is an epimorphism
between semigroups, where $\Omega_1,\ \Omega_2$ are bounded domains
in $\C^n,\ \C^m$ respectively, then $\varphi$ is an isomorphism. 
\end{theorem}

The author is grateful to A. Eremenko for his guidance and numerous 
suggestions concerning the paper. He also thanks L. Avramov, S. Bell 
and A. Gabrielov for valuable discussions and their interest in this work.

\section{Topology}\label{S:Top}

For a bounded domain $\Omega$ in $\C^n$ we denote by $C(\Omega)$ the 
subsemigroup of $E(\Omega)$ consisting of constant maps. An endomorphism 
$c_z$ is constant if it sends $\Omega$ to a point $z\in\Omega$. The subset 
$C(\Omega)\subset E(\Omega)$ can be described using only the semigroup 
structure as follows:

\begin{equation}\label{E:Co}
c\in C(\Omega) \ \text{iff} \ \ \forall (f\in E(\Omega)), \ \ (c\circ f=c).
\end{equation}

It is clear that we have a bijection between constant endomorphisms of 
$\Omega$ and points of this domain as a set: to each $z$ corresponds a 
unique $c_z\in C(\Omega)$ and vice versa, so we can identify the two. 
Under this identification, a subset of $\Omega$ corresponds to a 
subsemigroup of $C(\Omega)$.

Having defined points of a domain in terms of its semigroup structure of 
analytic endomorphisms, we can construct a map $\psi$ between $\Omega_1$ 
and $\Omega_2$ as follows

\begin{equation}\label{E:Dp}
\psi(z)= w \ \ \text{iff}\ \ \varphi c_z=c_w.
\end{equation}

So defined, $\psi$ satisfies (\ref{E:C}). Indeed, let $f\in E(\Omega_1),\, 
f(z)=\zeta
$. This is equivalent to 
\begin{equation}\label{E:Se}
f\circ c_z=c_{\zeta}. 
\end{equation}
Applying $\varphi$ to both sides of (\ref{E:Se}) we have 

\begin{equation}\label{E:3}
\varphi f\circ c_{\psi(z)}=c_{\psi(\zeta)}. 
\end{equation}
But (\ref{E:3}) is equivalent to $\varphi f(\psi(z))=\psi(\zeta)=\psi(f(z))$, 
which is (\ref{E:C}).

We describe the topology of a domain $\Omega$ using its injective 
endomorphisms. A map $f\in E(\Omega)$ is injective if and only if 
$$
\forall (c'\in C(\Omega))\ \forall (c''\in C(\Omega)),\ \ 
((f\circ c'= f\circ c'')\Rightarrow(c'= c'')).
$$
We denote the class of injective endomorphisms of $\Omega$ by 
$E_i(\Omega)$. For every $f\in E_i(\Omega),\ f_i(\Omega)$ is open 
\cite{sB48}.
The family
$
\{f(\Omega), \ \ f\in E_i(\Omega)\}
$
of subsets of $\Omega$
forms a base of topology, because every $z\in \Omega$ has a neighborhood
$f(\Omega)$, where $f(\zeta)=z + \lambda(\zeta-z)$, $f$ belongs to 
$E_i(\Omega)$ for every $\lambda$ such that $|\lambda|$ is small. 

To summarize, we described subsets of $\Omega$ and the topology on it 
using only the semigroup structure of $E(\Omega)$. Since this is so, the 
semigroup structure also defines the notions of an open set, closed set, 
compact set, closure of a set.

Now we can easily prove continuity of the map $\psi$ constructed above. 
Indeed, let $g(\Omega_2), \ g\in E_i(\Omega_2)$ be a set from the base of 
topology of $\Omega_2$. We take $f=\varphi^{-1}g$. Then  
$f\in E_i(\Omega_1)$ and $\psi^{-1}(g(\Omega_2))=f(\Omega_1)$, which 
proves that $\psi$ is continuous. Since $\varphi$ is an isomorphism, 
the same argument works to prove that $\psi^{-1}$ is also continuous,  
and thus $\psi$ is a homeomorphism.

Therefore the domains $\Omega_1, \ \Omega_2$ are homeomorphic, and hence 
\cite{wH48} they have the same dimension, i. e. $n = m$.

\section{Localization}\label{S:Loc}

We need the following lemma.
\begin{lemma}\label{L:E}
Suppose $H$ is a semigroup with identity, and $f$ an element of $H$ with 
the following two properties:

(i) $hf=fh$, for every $h$ in $H$;

(ii) $h_1f=h_2f$ implies $h_1=h_2$, for every $h_1$ and $h_2$ in $H$.

Then there exists a semigroup $S_f$ and a monomorphism $i:\ H\rightarrow 
S_f$, such that $i(f)$ is invertible in $S_f$ and commutes with all 
elements of 
$S_f$. Moreover, the semigroup $S_f$ satisfies the following universal 
property: for every semigroup $S_1$ with a monomorphism 
$i_1: \ H\rightarrow S_1$ such that $i_1(f)$ is invertible in $S_1$ 
and commutes with all elements of $S_1$, there exists a unique 
monomorphism $\hat i_1: \ S_f\rightarrow S_1$ such that 
$i_1=\hat i_1\circ i$.  
\end{lemma}

\begin{remark}
Uniqueness of $\hat i_1$ implies that the semigroup $S_f$ with the  
universal property is unique up to an isomorphism.
\end{remark}

\emph{Proof.}
We construct $S_f$ as follows. First we consider formal expressions of the 
form $hf^k$, where $h\in H$ and $k$ is an integer (may be positive, negative 
or 
zero). Then we define a multiplication on this set: 
$h_1f^{k_1}*h_2f^{k_2}=h_1h_2f^{k_1+k_2}$. Next we consider a relation 
on the set of formal expressions: $h_1f^{k_1}\sim h_2f^{k_2}$ if 
$k_1\leq k_2$ and $h_1=h_2f^{k_2-k_1}$ in $H$, or $k_2\leq k_1$ and 
$h_2=h_1f^{k_1-k_2}$ in $H$. It is easy to verify that this is an 
equivalence relation and it is compatible with the operation $*$; 
that is, $x\sim y, \ u\sim v$ implies $x*u\sim y*v$. 

Lastly, let $S_f$ be the set of equivalence classes with the binary 
operation induced by $*$. For $S_f$ to be a semigroup, we need to show 
that the binary operation $*$ is associative. Let 
$h_1f^{k_1}\sim h_1'f^{k_1'}$, $h_2f^{k_2}\sim h_2'f^{k_2'}$ and 
$h_3f^{k_3}\sim h_3'f^{k_3'}$. We need to show that 
$(h_1f^{k_1}*h_2f^{k_2})*h_3f^{k_3}\sim h_1'f^{k_1'}*
(h_2'f^{k_2'}*h_3'f^{k_3'})$. By the definition of the operation $*$, 
the last equivalence is the same as $h_1h_2h_3f^{k_1+k_2+k_3}
\sim h_1'h_2'h_3'f^{k_1'+k_2'+k_3'}$. Assuming that $k_1+k_2+k_3
\leq k_1'+k_2'+k_3'$, we have essentially one possibility to consider 
(the others are either similar or trivial): $k_1\leq k_1', \ k_2\leq 
k_2', \ k_3'\leq k_3$. In this case $h_1h_2h_3f^{k_3-k_3'}=
h_1'h_2'h_3'f^{k_1'-k_1+k_2'-k_2}$. Now we can use the cancellation 
property (ii) to get the desired equivalence.

The semigroup $H$ is embedded into $S_f$ via $i: \ h\mapsto [hf^0]$. 
The element $i(f)=[\text{id} f]$, where $\text{id}$ is the identity in 
$H$, is invertible in $S_f$ with the inverse $[\text{id} f^{-1}]$. 
Clearly, $[{\text{id}}f]$ commutes with all elements of $S_f$.

Now, suppose that $S_1$, $i_1:\ H\rightarrow S_1$ is a semigroup and a 
monomorphism, such that $i_1(f)$ is invertible in $S_1$ and commutes 
with all elements of $S_1$. Then we define 
$$
\hat i_1([hf^k])=i_1(h)(i_1(f))^k.
$$
This definition does not depend on a representative of $[hf^k]$. Indeed,
suppose $h_1f^{k_1}\sim h_2f^{k_2}$ and assume $k_1\leq k_2$. 
Then $h_1=h_2f^{k_2-k_1}$, and thus $i_1(h_1)=i_1(h_2)i_1(f)^{k_2-k_1}$. 
Hence $i_1(h_1)i_1(f)^{k_1}=i_1(h_2)i_1(f)^{k_2}$.

So defined, $\hat i_1$ is a homomorphism:

\begin{align}
\hat i_1 
&([h_1f^{k_1}][h_2f^{k_2}])=\hat i_1([h_1h_2f^{k_1+k_2}])\notag \\
&=i_1(h_1h_2)i_1(f)^{k_1+k_2}=i_1(h_1)i_1(h_2)i_1(f)^{k_1}i_1(f)^{k_2}\notag\\
&=i_1(h_1)i_1(f)^{k_1}i_1(h_2)i_1(f)^{k_2}=\hat i_1([h_1f^{k_1}])\hat i_1([h_2f^{k_2}]).\notag
\end{align}

The relation $\hat i_1\circ i= i_1$ holds, since $\hat i_1([hf^0])=i_1(h)$ 
for all $h\in H$.

Uniqueness of $\hat i_1$ is clear.
Lemma \ref{L:E} is proved.

\section{Extension of $\varphi$}\label{S:Ext}

Following \cite{aE93}, we say that for a bounded domain $\Omega$ an 
element $f\in E(\Omega)$ is {\it{good at}} $z\in \Omega$, denoted by 
$f\in G_z(\Omega)$, if 

\begin{enumerate}
\item $z$ is a unique fixed point of $f$;
\item $f(\Omega)$ has compact closure in $\Omega$; 
\item $f$ is injective in $\Omega$.
\end{enumerate}

Property 3 of a good element was already stated in terms of the semigroup 
structure of $\Omega$. Since the topology on $\Omega$ was described using 
only the semigroup structure, Property 2 can also be stated in these 
terms. Property 1 can be expressed in terms of the semigroup structure as
$$
(f\circ c_z=c_z)\wedge((f\circ c_{\zeta}=c_{\zeta})\Rightarrow 
(c_{\zeta}=c_z)).
$$

Since $f$ is an endomorphism of a domain, all eigenvalues $\lambda$ of its 
linear part at $z$ satisfy $|\lambda|\leq 1$ \cite{sK98}. Moreover, 
$|\lambda|< 1$ 
because the closure of $f(\Omega)$ is a compact set in $\Omega$. The 
injectivity of $f$ implies \cite{sB48} that it is biholomorphic onto 
$f(\Omega)$ and the Jacobian determinant of $f$ does not vanish at any 
point of $\Omega$. 

It is clear that for every $z\in \Omega$ a good element $f$ at $z$ exists. 
For example, we can take $f(\zeta)=z+\lambda (\zeta-z)$ with sufficiently 
small $|\lambda|$.

Consider a good element $f\in G_z(\Omega)$ and its commutant $H_f(\Omega)$ 
in $E(\Omega)$:
$$
H_f(\Omega)=\{h\in E(\Omega):\ \ hf=fh\}.
$$
Clearly $H_f(\Omega)$ is a subsemigroup of $E(\Omega)$. The element $f$, 
being good (hence injective), satisfies the cancellation property $(ii)$ 
of Lemma~\ref{L:E} in $H_f(\Omega)$. Thus, by Lemma~\ref{L:E}, 
we have the extension 
$S_f$ 
of $H_f(\Omega)$ in which $f$ is invertible and commutes with all elements 
of $S_f$. In the case of analytic endomorphisms we can embed $H_f(\Omega)$ 
into the subsemigroup of $A_z$, the semigroup of germs of analytic mappings 
at $z$ under composition, consisting of elements that commute with the germ 
of $f$ and containing the 
germ of $f^{-1}$. We use the universal property of Lemma~\ref{L:E} 
to conclude 
that $S_f$ is isomorphic to a subsemigroup of $A_z$. We identify $S_f$ with 
this semigroup, i. e. we consider elements of $S_f$ as germs of 
analytic mappings at $z$.
 
In proving that $\psi$ is (anti-) biholomorphic we need to show that it is 
so in a neighborhood of every point of $\Omega_1$. Since an (anti-) 
biholomorphic type of a domain is preserved by translations in $\C^n$, 
it is enough to show that $\psi$ is (anti-) biholomorphic in a neighborhood 
of $0\in \C^n$, assuming that $\Omega_1$ and $\Omega_2$ contain 0 and 
$\psi(0)=0$.

Let $\varphi:\ E(\Omega_1)\rightarrow E(\Omega_2)$ be an isomorphism of the 
semigroups, $f$ a good element, $f\in G_0(\Omega_1)$, and $H_f(\Omega_1)$ 
the commutant of $f$. Then clearly $H_g(\Omega_2)=\varphi(H_f(\Omega_1))$ 
is the commutant of $g=\varphi f$.
By Lemma~\ref{L:E}, we have the extensions $S_f, \ S_g$ of $H_f(\Omega_1)$ and 
$H_g(\Omega_2)$ respectively, and by the universal property of this lemma 
the isomorphism $\varphi$ extends to an isomorphism
$$
\Phi:\ \ S_f\rightarrow S_g.
$$

\section{System of projections and linearization}\label{S:Sys}

Let $\Omega$ be a bounded domain in $\C^n$. We say that a good  element 
$f\in G_0(\Omega)$ is {\it{very good at}} 0, and write $f\in VG_0(\Omega)$, 
if the corresponding semigroup $S_f\subset A_0$ constructed in 
Section~\ref{S:Ext}
contains a system of elements, which we call a system of projections, 
$\{p_i\}_{i=1}^{n}$ with the following properties:

(a) $\forall \ (i=1,\dots, n),\ \ (p_i\neq 0)$;

(b) $\forall \ (i=1,\dots, n),  \ \ (p_i^2=p_i)$;

(c) $\forall\ (i,\ j=1,\dots, n,\ i\neq j), \ \ (p_ip_j=0)$.

There does exist a very good element, since we can take $f$ to be a 
homothetic transformation at 0 with sufficiently small coefficient, $p_i$ a 
projection on the $i$'th coordinate of the standard coordinate system. 
Clearly, $p_i f = f p_i$ and there exists $k$ such that 
$p_if^k\in E(\Omega)$, and hence $p_i\in S_f$. From now on, we fix a very 
good element $f\in VG_0(\Omega)$, associated semigroups $H_f(\Omega),\ S_f$ 
and a system of projections $\{p_i\}$.

We introduce another subsemigroup of $E(\Omega)$:
$$
P_f(\Omega)=\{h\in G_0(\Omega)\cap H_f(\Omega),\ \ hp_i=p_ih\,\ \  i=1,\dots, 
n\},
$$
where the commutativity relations are in $S_f\subset A_0$. Notice that
$P_f(\Omega)\neq\emptyset$ since $f$ belongs to it. 

\begin{lemma}\label{L:Lin}
For every $h\in P_f(\Omega)$ there exists a biholomorphic germ $\theta_h$ at 
$0\in \C^n$ such that $\theta_h h=\Lambda\theta_h$, where $\Lambda=
{\rm{diag}}(\lambda_1,\dots, \lambda_n)$ is an invertible diagonal matrix 
which is similar to $dh(0)$ in $GL(n,\C)$.
\end{lemma}

\emph{Proof.}
The relations $p_i\neq0,\ p_i^2=p_i,\ p_ip_j=0, \ i\neq j$, imply that for 
$P_i=dp_i(0)$, the linear part of $p_i$ at 0, we have 
$P_i\neq 0, \ P_i^2=P_i, \ P_iP_j=0, \ i\neq j$. Since the matrices $P_i$ 
commute, there exists \cite{kH71} a matrix $A\in GL(n, \C)$ such that 
$P_i'=AP_iA^{-1}=\Delta_i=\text{diag}(0,\dots, 1,\dots, 0)$, where the only 
non-zero entry appears in the $i$'th place.

Since $p_i^2=p_i,\ i=1,\dots, n$, we can use the argument given in 
\cite{sK98} to linearize $p_i$, i. e. there exists a biholomorphic 
germ $\xi_i$ at 0 such that $\xi_ip_i=P_i\xi_i,\ d\xi_i(0)=\text{id}, \ 
i=1,\dots, n$. The map $\xi_i$ is constructed in \cite{sK98} as follows:
$$
\xi_i=\text{id}+(2P_i-\text{id})(p_i-P_i), \ \ i=1,\dots, n.
$$ 
If we take $\xi_i'=A\xi_i$, we have $\xi_i'p_i=P_i'\xi_i'$.
For simplicity of notations, we assume that $\xi_i$ itself conjugates $p_i$ 
to a diagonal matrix, that is, $P_i=P_i'$ (in this case $P_i$ is not 
necessarily $dp_i(0)$, but rather $Adp_i(0)A^{-1};\ d\xi_i(0)=A$).
For every $i=1,\dots, n$ we have $h_iP_i=P_ih_i$, where 
$h_i=\xi_ih\xi_i^{-1}$. Let $H_i=dh_i(0)$. Then $H_iP_i=P_iH_i$, and hence 
in the $i$'th row and the $i$'th column the matrix $H_i$ has only one 
non-zero entry, $\lambda_i$, which is located at their intersection. 
Thus $\lambda_i$ has to be an eigenvalue of $H_i$, and hence of the linear 
part of $h$. In particular, $0<|\lambda_i|<1$.

Let $I_i:\ \C\rightarrow\C^n$ be the embedding $z\mapsto 
(0,\dots, z, \dots, 0)$, where the only non-zero entry is $z$, which is 
in the $i$'th place; and $\pi_i:\ \C^n\rightarrow \C$, a projection 
$(z_1,\dots, z_n) \mapsto z_i$, corresponding to the $i$'th axis. For 
every $i=1,\dots, n$, the map $\pi_ih_iI_i$ sends a neighborhood of 0 in 
$\C$ into $\C$, and its derivative at 0, $ \lambda_i$, is an 
eigenvalue of $h$. Hence (\cite{lC93}, p. 31) $\pi_ih_iI_i$ is linearized 
by the 
unique solution $\eta_{h,i}$ of the Schr\"oder equation 
\begin{equation}\label{E:4}
\eta(\pi_ih_iI_i)=\lambda_i\eta, \ \ \eta(0)=0,\ \eta'(0)=1. 
\end{equation}
Since $P_iI_i=I_i, \ \pi_iP_iI_i=\text{id}_{\C}$, we can rewrite (\ref{E:4}) 
as 
$$
\eta_{h,i}\pi_ih_iP_iI_i=\lambda_i\eta_{h,i}\pi_iP_iI_i,\ \  {\text{or}}\ \ 
\eta_{h,i}\pi_ih_iP_i=\lambda_i\eta_{h,i}\pi_iP_i.
$$
But $h_iP_i=P_ih_i$, and so 
\begin{equation}\label{E:5}
\eta_{h,i}\pi_iP_ih_i=\lambda_i\eta_{h,i}\pi_iP_i. 
\end{equation}
The equation (\ref{E:5}), in its turn, is equivalent to 
\begin{equation}\label{E:6}
\eta_{h,i}\pi_i\xi_ip_ih=\lambda_i\eta_{h,i}\pi_i\xi_ip_i. 
\end{equation}
We denote  
\begin{equation}\label{E:7}
\theta_{h,i}=\eta_{h,i}\pi_i\xi_ip_i, 
\end{equation}
a map from a neighborhood of $0\in\C^n$ into $\C$.
Then (\ref{E:6}) becomes $\theta_{h,i}h=\lambda_i\theta_{h,i}$. Now we define
$$
\theta_h=(\theta_{h,1},\dots, \theta_{h,n}),
$$ 
which is a germ of an analytic map at 0. This germ linearizes $h$:
$$
\theta_hh=(\theta_{h,1}h,\dots,\theta_{h,n}h)=(\lambda_1\theta_{h,1},\dots, \lambda_n\theta_{h,n})=\Lambda \theta_h,
$$
where $\Lambda=\text{diag}(\lambda_1,\dots, \lambda_n)$ is an invertible diagonal matrix, which has eigenvalues of $dh(0)$ on its diagonal.

The germ $\theta_h$ is biholomorphic. Indeed, 
$$
\theta_{h,i}=\eta_{h,i}\pi_i\xi_ip_i=\eta_{h,i}\pi_iP_i\xi_i,\ \ i=1,\dots, n.
$$ 
Using the chain rule, we see that $d\theta_h(0)=A$, where $A$ is an invertible diagonal matrix that diagonalizes $P_i$. We conclude that $\theta_h$ is biholomorphic. Lemma~\ref{L:Lin} is proved.

\section{Simultaneous linearization }\label{S:Sim}

Using Lemma~\ref{L:Lin}, we can linearize elements of $P_f(\Omega)$. Namely, for every 
$h\in P_f(\Omega)$ there exists $\theta_h$ (constructed in 
Section~\ref{S:Sys}), such that $\theta_hh=\Lambda_h\theta_h$, where 
$\Lambda_h$ is an invertible diagonal matrix. In particular, we can 
linearize $f$:
$$
\theta_ff=\Lambda_f\theta_f,
$$
where the germ $\theta_f$ is biholomorphic at 0, and $\Lambda_f$ is an 
invertible diagonal matrix. 

\begin{lemma}\label{L:Slin}
For every $h\in P_f(\Omega)$ we have $\theta_h=\theta_f$.
\end{lemma}

\emph{Proof.}
Let us consider the germ 
\begin{equation}\label{E:8}
\theta=\Lambda_f^{-1}\theta_hf, 
\end{equation}
which is clearly biholomorphic. We have
$$
\theta h=\Lambda_f^{-1}\theta_hfh=\Lambda_f^{-1}\theta_hhf=\Lambda_f^{-1}
\Lambda_h\theta_hf=\Lambda_h\Lambda_f^{-1}\theta_hf=\Lambda_h\theta.
$$ 
Using (\ref{E:8}), we write the equation $\theta h=\Lambda_h\theta$ in the 
coordinate form:
$$
(1/{\lambda_{f,i}})\theta_{h,i}fh=({\lambda_{h,i}}/{\lambda_{f,i}})
\theta_{h,i}f, \ \ i=1,\dots, n.
$$
By (\ref{E:7}) and the definition of $\xi_i$,
$$
(1/{\lambda_{f,i}})\eta_{h,i}\pi_iP_if_ih_i=({\lambda_{h,i}}/{\lambda_{f,i}})
\eta_{h,i}\pi_i P_if_i,\ \ i=1,\dots, n,
$$
where $f_i=\xi_if\xi_i^{-1}$.
Using the commutativity relations $f_iP_i=P_if_i$,  $h_iP_i=P_ih_i$, which hold 
since $\{p_i\}\subset S_f, \ h\in P_f(\Omega)$, we get
\begin{align}
&(1/{\lambda_{f,i}})\eta_{h,i}\pi_if_ih_iP_i=({\lambda_{h,i}}/{\lambda_{f,i}})
\eta_{h,i}\pi_if_iP_i,\ \ \text{or} \notag\\
&(1/{\lambda_{f,i}})\eta_{h,i}\pi_if_ih_iI_i=({\lambda_{h,i}}/{\lambda_{f,i}})
\eta_{h,i}\pi_if_iI_i, \ \ i=1,\dots, n.\notag
\end{align}
This is the same as 
$$
((1/{\lambda_{f,i}})\eta_{h,i}\pi_if_iI_i)(\pi_ih_iI_i)=\lambda_{h,i}
((1/{\lambda_{f,i}})\eta_{h,i}\pi_if_iI_i), \ \ i=1,\dots, n,
$$
since $h_i$ locally preserves the $i$'th coordinate axis ($h_iP_i=P_ih_i$).
It is easily seen that 
\begin{align}
&((1/{\lambda_{f,i}})\eta_{h,i}\pi_if_iI_i)(0)=0, \notag \\
&((1/{\lambda_{f,i}})\eta_{h,i}\pi_if_iI_i)'(0)=1.\notag
\end{align}
A normalized solution to a Schr\"oder equation is unique, though; thus we have
$$
\eta_{h,i}(\pi_if_iI_i)=\lambda_{f,i}\eta_{h,i}, \ \ \eta_{h,i}(0)=0,
\ \eta_{h,i}'(0)=1.
$$
Using the uniqueness argument again, we obtain $\eta_{h,i}=\eta_{f,i}$, and 
hence $\theta_h=\theta_f$. The lemma is proved.

According to Lemma~\ref{L:Slin}, the single biholomorphic germ $\theta_f$ 
conjugates the subsemigroup $P_f(\Omega)$ to some subsemigroup $D_f$ of invertible 
diagonal matrices in $D_n$, the set of all $n\times n$ diagonal matrices 
with entries in $\C$. We show that $D_f$ contains all invertible 
diagonal matrices with sufficiently small entries. To do this, first we 
extend $\theta_f$ to an analytic map on the whole domain $\Omega$ using 
the formula
$$
\theta_f=\Lambda_f^{-l}\theta_ff^l,
$$
where $l$ is chosen so large that $\text{Cl}\{f^l(\Omega)\}$ is contained in 
a neighborhood of 0 where $\theta_f$ is originally defined and biholomorphic; 
the symbol Cl denotes closure. From the procedure of extending $\theta_f$ to 
$\Omega$ we see that it is one-to-one and bounded in the whole domain.

Now, let $\Lambda=\text{diag}(\lambda_1,\dots,\lambda_n)$ be a matrix such 
that $\text{Cl}\{\Lambda\theta_f(\Omega)\}\subset W$, where $W$ is a 
neighborhood of $0\in\C^n$ for which $\text{Cl}\{\theta_f^{-1}W\}
\subset\Omega$. Such a matrix $\Lambda$ exists since $\theta_f$ is bounded 
in $\Omega$. Consider $h=\theta_f^{-1}\Lambda\theta_f$, which belongs to 
$G_0(\Omega)$. The map $h$ commutes with $f$ and all $p_i$'s. Indeed, using 
the formula $\theta_ff\theta_f^{-1}=\Lambda_f$, we conclude that $hf=fh$ is 
equivalent to $\Lambda\Lambda_f=\Lambda_f \Lambda$, which is a true relation 
since both matrices $\Lambda$ and $\Lambda_f$ are diagonal. The relations 
$hp_i=p_ih, \ i=1,\dots, n$, are verified similarly, using the formula 
$\theta_fp_i\theta_f^{-1}=P_i$, which follows from the definition of 
$\theta_f$.

\section{Solving a matrix equation}\label{S:Sol}

We proved that for an element $f\in VG_0(\Omega)$ there exists a 
biholomorphic germ $\theta_f$ conjugating the semigroup $P_f(\Omega)$ 
to a subsemigroup $D_f\subset D_n$, which contains all invertible diagonal 
matrices with sufficiently small entries.

Let $f\in VG_0(\Omega_1)$, and $g=\varphi f$. Then $g\in VG_0(\Omega_2)$, 
and there is an isomorphism
$$
\Phi:\ \ S_f\rightarrow S_g.
$$
For the mappings $f$ and $g$ we have 
$$
\theta_ff=\Lambda_f\theta_f,\ \ \theta_gg=M_g\theta_g,
$$
where $\Lambda_f, \ M_g$ are invertible diagonal matrices. 

Let us consider the germ $L=\theta_g\psi \theta_f^{-1}$. This germ conjugates 
the semigroups $D_f, \ D_g$:
\begin{align}
L\Lambda L^{-1}
&=\theta_g\psi\theta_f^{-1}\Lambda\theta_f\psi^{-1}
\theta_g^{-1} \notag\\ 
&=\theta_g\psi h\psi^{-1}\theta_g^{-1}=\theta_g j 
\theta_g^{-1}=M,\notag
\end{align}
where $h\in P_f, \ \theta_fh=\Lambda\theta_f;\ j=\varphi h,\ \theta_gj=M
\theta_g$. 

Define $R(\Lambda)=L\Lambda L^{-1}$. Then $R:\ D_f\rightarrow D_g,$
$$
R(\Lambda_1\Lambda_2)=R(\Lambda_1)R(\Lambda_2), \ \ \Lambda_1, \ \Lambda_2
\in D_f.
$$

In what follows, we will identify $D_n$ with the multiplicative semigroup 
$\C^n$ ($D_n\cong \C^n$) in the obvious way and consider a topology 
on $D_n$ induced by the standard topology on $\C^n$.

We are going to extend $R$ to an isomorphism of $D_n$. First, 
we denote by $\overline{D}_f,\ \overline{D}_g$ the closures of $D_f,\ D_g$ 
in $D_n$, and for $\Lambda\in \overline{D}_f$ we set 
$$
R(\Lambda)=\lim R(\Lambda_k),\ \ \Lambda_k\rightarrow\Lambda,\ \Lambda_k
\in D_f.
$$
This limit exists and does not depend on the sequence $\{\Lambda_k\}$, which 
follows from the fact that $\psi^{\pm1}, \theta_f^{\pm1}, \theta_g^{\pm1}$ 
are continuous. The map $R$ is an isomorphism of topological semigroups 
$\overline{D}_f$ and $\overline{D}_g$ (the inverse of $R$ has a similar 
representation).

Next, we extend the map $R$ to $D_n$ as
$$
R(\Gamma)=R(\Gamma\Lambda)R(\Lambda)^{-1}, \ \ \Gamma\in D_n,
$$
where $\Lambda\in D_f$ is chosen so that $\Gamma\Lambda\in\overline{D}_f$. 
This definition does not depend on the choice of $\Lambda$. Indeed, since 
all matrices in question are diagonal (hence commute), the relation 
$R(\Gamma\Lambda_1)R(\Lambda_1)^{-1}=R(\Gamma\Lambda_2)R(\Lambda_2)^{-1}$ 
is equivalent to $R(\Gamma\Lambda_1)R(\Lambda_2)=R(\Gamma\Lambda_2)
R(\Lambda_1)$, which holds. 

The extended map $R$ is clearly an isomorphism of $D_n$ onto itself. Thus 
we have
\begin{equation}\label{E:9}
R(\Lambda'\Lambda'')=R(\Lambda')R(\Lambda''), \ \ \Lambda', \ \Lambda''\in 
D_n. 
\end{equation}

Injectivity of $R$ and (\ref{E:9}) imply that $R(\Delta_i)=\Delta_j$ for 
all $i$, 
where $j=j(i)$ depends on $i$; $j(i)$ is a permutation on $\{1,\dots, n\}$ 
(we recall that $\Delta_i=\text{diag}(0,\dots,1,\dots,0)$). This is because 
$\{\Delta_i\}_{i=1}^n$ is the only system in $D_n$ with the following 
relations: $\Delta_i\neq 0,\ \Delta_i^2=\Delta_i,\ \Delta_i\Delta_j=0, 
\ i\neq j$. 

Since all matrices $\Lambda$ and their images $R(\Lambda)$ are diagonal, we 
can consider the matrix equation (\ref{E:9}) as $n$ scalar equations:
\begin{equation}\label{E:10}
r_j(\lambda_1'\lambda_1'',\dots, \lambda_n'\lambda_n'')=r_j(\lambda_1',
\dots, \lambda_n')r_j(\lambda_1'',\dots, \lambda_n''), \ \ j=1,\dots, n,
\end{equation}
where $r_j$ are components of $R$.
If we rewrite the equation $R(\Delta_i\Lambda)=\Delta_jR(\Lambda)$ in the 
coordinate form, we see that 
$$r_j(\lambda_1,\dots, \lambda_n)=r_j(0,\dots, 
\lambda_i,\dots, 0)=q_j(\lambda_i); 
$$
that is, each $r_j$ depends on only one 
of the $\lambda_i$'s.
For each $j$ the corresponding equation in (\ref{E:10}) in terms of the 
$q_j$'s 
becomes
$$
q_j(\lambda_i'\lambda_i'')=q_j(\lambda_i')q_j(\lambda_i'').
$$
This equation has (\cite{aE93}, p. 130) either the constant solution 
$q_j(\lambda_i)=1$, or 
$$
q_j(\lambda_i)=\lambda_i^{\alpha_{ij}}\overline{\lambda}_i^{\beta_{ij}}, 
\ \ \alpha_{ij},\  \beta_{ij}\in \C, \ \ \alpha_{ij}-\beta_{ij}=\pm1.
$$

Going back to the function $L$, we have
\begin{equation}\label{E:11}
\begin{aligned}
L\text{diag}(\lambda_1,\dots, \lambda_n)=\text{diag}
(\lambda_{i(1)}^{\alpha_1}\overline{\lambda}_{i(1)}^{\beta_1},
\dots, \lambda_{i(n)}^{\alpha_n}\overline{\lambda}_{i(n)}^{\beta_n})L, \\
 \alpha_i-\beta_i=\pm1, \ \ i=1,\dots, n,\notag
\end{aligned}
\end{equation}
where $i(j)$ is the inverse permutation to $j(i)$.

Let us choose and fix $(\mu_1,\dots, \mu_n)$ such that $(1/\mu_1, 
\dots, 1/\mu_n)$ belongs to a neighborhood $W_0$ of $0\in\C^n$ 
where $L$ is defined, and let $W_1$ be a neighborhood of $0\in\C^n$ 
such that $(\mu_1z_1,\dots, \mu_nz_n)\in W_0$, whenever 
$(z_1,\dots, z_n)\in W_1$.
Then from (\ref{E:11}) we have
\begin{align} 
L&(z_1,\dots, z_n)
=L\text{diag}(\mu_1 z_1,\dots, \mu_n z_n)(1/\mu_1,\dots, 1/\mu_n) \notag\\
&=\text{diag}((\mu_{i(1)}z_{i(1)})^{\alpha_1}(\overline{\mu_{i(1)} 
z_{i(1)}})^{\beta_1},\dots, (\mu_{i(n)}z_{i(n)})^{\alpha_n}
(\overline{\mu_{i(n)} z_{i(n)}})^{\beta_n}) \notag\\
&\times L (1/\mu_1,\dots, 1/\mu_n)=B (z_1^{\alpha_1}
\overline{z}_1^{\beta_1}, \dots, z_n^{\alpha_n}\overline{z}_n^{\beta_n}), 
\notag 
\end{align}
where $B$ is a constant matrix. The last formula is the explicit 
expression for $L$.

\section{Proving that $\psi$ is (anti-) biholomorphic}\label{S:Pr}

To prove that $\psi$ is (anti-) biholomorphic is the same as to prove that 
$L$ is (anti-) biholomorphic, because the relation $L=\theta_g\circ\psi\circ
\theta_f^{-1}$ holds. We showed that 
\begin{equation}\label{E:12}
L(z_1,\dots, z_n)=B (z_1^{\alpha_1}\overline{z}_1^{\beta_1}, \dots, 
z_n^{\alpha_n}\overline{z}_n^{\beta_n}), \ \  \alpha_i-\beta_i=\pm1, 
\ \ i=1,\dots, n 
\end{equation}
in a neighborhood $W_1$ of $0$. From the representation (\ref{E:12}) we see 
that $L$ is $\R$-differentiable and non-degenerate in $W_1\setminus 
\cup_{k=1}^n
\{(z_1, \dots, z_n):\ z_k=0\}$. Since this is true for every point in the 
domain $\Omega_1$, the map $\psi$ is $\R$-differentiable and 
non-degenerate everywhere, with the possible exception of an analytic set. 
Let us remove this set from $\Omega_1$, as well as its image under $\psi$ 
from $\Omega_2$. We call the domains obtained in this way $\Omega',
\ \Omega''$. Now the map $\psi:\ \Omega'\rightarrow\Omega''$ is 
$\R$-differentiable and non-degenerate everywhere. It is clear that 
if we prove that $\psi$ is (anti-) biholomorphic between $\Omega', 
\ \Omega''$, then it is (anti-) biholomorphic between $\Omega_1,\  \Omega_2$ 
due to a standard continuation argument \cite{sK82}. So we can think that 
$\psi$ is $\R$-differentiable and non-degenerate in $\Omega_1$ itself. 
The map $L$ thus has to be $\R$-differentiable and non-degenerate at 0. 
However, this is the case if and only if $\alpha_i+\beta_i=1, 
\ i=1,\dots, n$. Together with the equation $\alpha_i-\beta_i=\pm1$ it gives 
us that either $\alpha_i=1,\ \beta_i=0$, or $\alpha_i=0, \ \beta_i=1$. 

It remains to show that either $\alpha_i=1$ and $\beta_i=0$, or $\alpha_i=0$ 
and $\beta_i=1$, simultaneously for all $i$. Suppose, by way of contradiction, 
that we have $L(z_1,\dots, z_n)=B (\dots, z_i,\dots, \overline{z}_j,\dots)$. 
Then 
$$
L^{-1}(w_1,\dots, w_n)=(\dots, l_i(w_1,\dots, w_n),\dots, 
l_j(\overline{w}_1,\dots, \overline{w}_n),\dots),
$$
where $l_i,\ l_j$ are linear analytic functions. Let us look at an 
endomorphism $f_0$ of $\Omega_1$ in the form 
$$
f_0=\theta_f^{-1}\lambda(\dots, \theta_{f,i}\theta_{f,j},\dots, 
\theta_{f,j},\dots)\theta_f,
$$
where $\theta_{f,i}\theta_{f,j}$ is in the $i$'th place and $\theta_{f,j}$ 
in the $j$'th; $|\lambda|$ is sufficiently small.
Using (\ref{E:C}) and the definition of $L$, we have
$$
\theta_g\varphi f_0\theta_g^{-1}=\theta_g\psi f_0\psi^{-1}\theta_g^{-1}=
L\theta_f f_0 \theta_f^{-1}L^{-1}.
$$
So, 
\begin{align}
\theta_g\varphi &f_0\theta_g^{-1}(w_1,\dots, w_n) \notag\\
&= B' (\dots, l_i(w_1,\dots, w_n) l_j(\overline{w}_1,\dots, \overline{w}_n),
\dots, \overline{l}_j(w_1,\dots, w_n),\dots)\notag
\end{align}
for some constant matrix $B'$. This map, and hence $\varphi f_0$, is not 
analytic, though, in a neighborhood of 0, which is a contradiction. Thus $L$, 
and hence $\psi$, is either analytic or antianalytic in a neighborhood of 0. 

Theorem~\ref{T:Mt} is proved completely.

\section{Proof of Theorem~\ref{T:Gt}}\label{S:Pr2}

Since $\varphi$ is an epimorphism, it takes constant endomorphisms of
$\Omega_1$ to constant endomorphisms of $\Omega_2$, which follows
from~(\ref{E:Co}). Thus we can define a map $\psi:\ \Omega_1\to\Omega_2$
as in~(\ref{E:Dp}). Following the same steps as in verifying~(\ref{E:C}),
we obtain
\begin{equation}\label{E:Sc}
\varphi f\circ\psi = \psi\circ f,\ \ {\text{for all }} f\in E(\Omega_1).
\end{equation}
We will show that~(\ref{E:Sc}) implies bijectivity of $\psi$. The map 
$\psi$ is surjective. Indeed, let $w\in\Omega_2$ and $c_w$ be the 
corresponding constant endomorphism. Since $\varphi$ is an 
epimorphism, there exists $f\in E(\Omega_1)$, such that 
$\varphi f= c_w$. If we plug this $f$ into~(\ref{E:Sc}), we get
$$
\psi f(z) = w
$$
for all $z\in \Omega_1$. Thus $\psi$ is surjective.

To prove that $\psi$ is injective, we show that for every 
$w\in \Omega_2$, the full preimage of $w$ under $\psi$, $\psi^{-1}(w)$,
consists of one point. 

Assume for contradiction that $S_w=\psi^{-1}(w)$ consists of more than
one point for some $w\in\Omega_2$. The set $S_w$ cannot be all of $\Omega_1$, 
since $\psi$ is surjective. For $z_0\in\partial S_w\cap\Omega_1$ we
can find $z_1\in S_w$ and $\zeta\notin S_w$ which are arbitrarily close
to $z_0$. Let $z_2$ be a fixed point of $S_w$ different from $z_1$.
Consider a homothetic transformation $h$ such that $h(z_1)=z_1,\ h(z_2)=
\zeta$. Since the domain $\Omega_1$ is bounded, we can choose points
$z_1$ and $\zeta$ sufficiently close to each other so that $h$ belongs to
$E(\Omega_1)$. Applying~(\ref{E:Sc}) to $h$ we obtain
\begin{align}
&\varphi h(w)=\varphi h\circ \psi(z_1)=\psi\circ h(z_1)=\psi(z_1)=w; \notag\\  
&\varphi h(w)=\varphi h\circ \psi(z_2)=\psi\circ h(z_2)=\psi(\zeta)\neq w. 
\notag
\end{align}
The contradiction shows injectivity of $\psi$. Thus we have proved that
$\psi$ is bijective.

According to~(\ref{E:Sc}) we have
\begin{equation}\label{E:C1}
\varphi f=\psi\circ f\circ \psi^{-1}, \ \ \text{for all}\ f\in E(\Omega_1),
\notag
\end{equation} 
which implies that $\varphi$ is an isomorphism. 

Theorem~\ref{T:Gt} is proved.

\end{document}